\documentclass[12pt]{article}

\usepackage[english]{babel}
\usepackage{amsfonts}
\usepackage{amsmath}
\usepackage{amssymb}
\usepackage{mathrsfs}
\usepackage{latexsym}
\usepackage{multicol}
\usepackage{fancybox}
\usepackage{dsfont}
\usepackage{color}
\usepackage{hyperref}
\usepackage{subfigure} 
\usepackage{graphicx} 
\usepackage{caption}

\newtheorem{thm}{Theorem}[section]

\newtheorem{propn}[thm]{Proposition}

\newtheorem{rem}[thm]{Remark}

\newtheorem{defn}[thm]{Definition}

\def\Ai{\textrm{Ai}}
\def\Bi{\textrm{Bi}}
\def\Hi{\textrm{Hi}}
\def\Gi{\textrm{Gi}}

\begin{document}

\begin{center}
\Large{\bf{Notes on the zeros of the solutions of the non-homogeneous Airy's equation}}
\end{center}

\begin{center}
{ {Federico Zullo \\ DICATAM, University of Brescia, \\ Via Valotti 9, 25133, Brescia (IT) federico.zullo@unibs.it
 }}
\end{center}

\medskip
\medskip

\begin{abstract}
\noindent We present some observations on the distribution of the zeros of solutions of the nonhomogeneous Airy's equation. 
We show the existence of a principal family of solutions, with simple zeros, and particular solutions, 
characterized by a double zero in a given position of the complex plane.  A recursion, 
describing the distribution of the zeros, is introduced and the limits of its applicability are discussed. The results can be considered
a generalization of previous works on the distribution of the zeros for the solutions of the 
corresponding homogeneous equation. 
\end{abstract}

\bigskip\bigskip

\noindent


%\body

%\tableofcontents

\section{Introduction}\label{ra_sec1}
It is well known that the knowledge of the distribution of the zeros in the complex plane for transcendental entire functions of finite order gives the possibility to get a representation of such functions through the Weierstrass-Hadamard factorization theorem. In general, for entire functions defined by the solutions of suitable differential or difference equations, the problem of the localization of the corresponding zeros may be very demanding. For second order homogeneous differential equations the problem has been investigated under different perspectives and many results are collected in the literature \cite{H} \cite{Laine}. For entire solutions of  nonhomogeneous differential equations, there are many withstanding problems about the distribution and localization of the corresponding  solutions. In previous works \cite{Z1} \cite{Z2}, we examined the distribution of the zeros of the general solution of the homogeneous Airy differential equation through the introduction of two suitable parameters into the equation. A limiting process on a parametric recursion has been shown to give the location of the zeros and the dependence of the zeros by varying the parameters has been investigated. In this work we consider the solutions of a nonhomogeneous Airy equation obtained by adding a constant term to the homogeneous equation. As we will see, this simple extension considerably changes the structure of the solutions regarding the distribution of the zeros. 

As a last remark, we would like to underline that the nonhomogeneous Airy equation possesses several applications in mathematical physics \cite{AM}. For example, it has a strict relation with the second member of the Burger's hierarchy \cite{KS}:
\begin{equation}\label{burg}
\psi_t+(\psi_{xx}+\psi^3+3\psi\psi_x)_x=0.
\end{equation}
Under the self-similarity transformation
\begin{equation}
\psi=\frac{1}{\eta(3t)^{1/3}}f(z), \quad z\doteq \frac{x}{\eta(3t)^{1/3}}-\frac{b}{\eta^3}, \quad \eta^3=a
\end{equation}
equation (\ref{burg}) becomes
\begin{equation}
\frac{d^3f}{dz^3} +3f\frac{d^2f}{dz^2}+3f^2\frac{df}{dz}=(az+b)\frac{df}{dz}+af.
\end{equation}
Integrating once one gets
\begin{equation}\label{eqks}
\frac{d^2f}{dz^2}+3\frac{df}{dz}f+f^3=k+(az+b)f(z),
\end{equation}
where $k$ is the integration constant.  The previous equation will be considered in the next sections (see in particular equation (\ref{equ})). Some of the solutions of (\ref{eqks}) have been considered in  \cite{KS} for the description of liquids with gas bubbles.

\section{The nonhomogeneous equation and a particular family of solutions}\label{sec1}
The nonhomogeneous Airy differential equations
\begin{equation}\label{eq}
\frac{d^2y(z)}{dz^2}=zy(z)\pm\frac{1}{\pi},
\end{equation}
and their particular solutions, denoted by $\Hi(z)$ and $\Gi(z)$ (the Scorer functions), have been considered, among others, by  Scorer \cite{S}, Olver \cite{Olver}, MacLeod \cite{AM}, Gil et al. \cite{GST} \cite{GST1}. In particular, the functions $\Hi(z)$ possesses the integral representation
\begin{equation}%\label{eq}
\Hi(z)=\frac{1}{\pi}\int_0^{\infty}\exp\left(zt-\frac{t^3}{3}\right)dt,
\end{equation}
and is a solution of (\ref{eq}) with the $+$ sign, whereas $\Gi(z)$ has the representation
\begin{equation}%\label{eq}
\Gi(z)=\frac{1}{\pi}\int_0^{\infty}\sin\left(zt+\frac{t^3}{3}\right)dt,
\end{equation}
and is a solution of (\ref{eq}) with the $-$ sign. The corresponding initial conditions are given by \cite{NIST}
\begin{equation}%\label{eq}
\Gi(0)=\frac{1}{2}\Hi(0)=\frac{1}{3^{7/6}\Gamma(2/3)},\quad \Gi'(0)=\frac{1}{2}\Hi'(0)=\frac{1}{3^{5/6}\Gamma(1/3)}
\end{equation}
The two functions are related by the equation $\Gi(z)+\Hi(z)=\Bi(z)$, where $\Bi(z)$ is one of the two fundamental solutions of the homogeneous Airy equation \cite{NIST}. Both $\Gi$ and $\Hi$ are entire functions with an infinite number of zeros in the complex plane: many properties of the zeros of these two functions and their derivatives can be found  in \cite{GST1}. Also, $\Hi$ decays at infinity on the negative real axis, whereas $\Gi$ decays on the positive real axis and is always less than 1 in magnitude on the real axis.  

To not keep the $\pm$ signs, let us insert a constant $c$ instead of $\pm \frac{1}{\pi}$ in (\ref{eq}):
\begin{equation}\label{eq1}
\frac{d^2y(z)}{dz^2}=zy(z)+c,
\end{equation}
It is clear that with a rescaling of the dependent variable it is possible to fix the value of the constant $c$ to be any suitable value. The general solution of equation (\ref{eq1}) can be written as
\begin{equation}\label{gs}
y(z)=C_1 \Ai(z)+C_2\Bi(z)+c\pi\Hi(z)
\end{equation}
where $C_1$ and $C_2$ are two arbitrary constants. $\Ai$, $\Bi$ and $\Hi$ are entire with order equal to $3/2$, so equation (\ref{eq1}) possesses entire solutions with order of growth equal to $\frac{3}{2}$: since this number is not an integer, all solutions of (\ref{eq1}) possesses an infinite number of zeros \cite{T}. These zeros are movable, in the sense that if the initial conditions change, the positions of the zeros change. The position of one zero can be considered a constant of integration. A second constant is the value of the derivative of the function at this zero. This is true also for the homogeneous equation corresponding to equation (\ref{eq}) (see \cite{Z1} \cite{Z2}), but in this case there is a significant difference. It is known that, further to the trivial solution, entire solutions of linear, second order, homogeneous differential equations 
\begin{equation}\label{lin}
\frac{d^2y(z)}{dz^2}=A(z)\frac{dy(z)}{dz}+B(z)y(z),
\end{equation}
where $A(z)$ and $B(z)$ are entire functions, have all their zeros simple. This result comes directly from the local analysis of equation (\ref{lin}) around any zero $z_0$. In the nonhomogeneous case one has the possibility to have double zeros.  Given any $p \in \mathbb{C}$, equation (\ref{eq1}) possesses always a solution with a double zero in $p$. Indeed it holds the following
\begin{propn}\label{p1}
Given $p \in \mathbb{C}$, equation (\ref{eq1}) possesses just one solution $\tau_d(z,p)$, proportional to $c$, with a double zero in $z=p$. This solution is defined by the series
\begin{equation}\label{dz}
\tau_d(z,p)=c\sum_{n=2}e_n(z-p)^n, \quad e_2=\frac{1}{2}, \; e_3=0, \; e_4=\frac{p}{24},\; e_{n+2}=\frac{pe_n+e_{n-1}}{(n+1)(n+2)}
\end{equation}
For $p$ real, $\tau_d(z,p)$ can also be represented as
\begin{equation}
\tau_d(z,p)=c\pi \left(\textnormal{Bi}(z)\int_{p}^z \textnormal{Ai}(x)dx-\textnormal{Ai}(z)\int_{p}^z\textnormal{Bi}(x)dx\right)
\end{equation}
\end{propn}
The second part of this proposition is also given in \cite{GST1} , the first part can be obtained by considering the Taylor series of $y(z)$ with a double zero in $z=p$ around this point. As we will see, the solutions $y(z)$ possessing a zero in $z=0$ (double or simple) will play a central role in the rest of the paper. In the case $p=0$ the recursion in (\ref{dz}) can be solved and the corresponding series is a generalized hypergeometric series. More explicitly one has
\begin{equation}\label{dz0}
\tau_d(z,0)=c\sum_{n=0}\frac{3^n n!}{(3n+2)!}z^{3n+2} = \frac{c}{2}z^2 {}_1F_2\left(1; \frac{4}{3},\frac{5}{3};\frac{z^3}{9}\right).
\end{equation}
The hypergeometric function $_1F_2\left(1; \frac{4}{3},\frac{5}{3};\frac{z^3}{9}\right)$ is actually related to the Lommel functions $s_{\mu,\nu}(z)$. The Lommel functions are particular solutions of the nonhomogeneous Bessel equations
\begin{equation}\label{nhB}
z^2\frac{d^2y(z)}{dz^2}+z\frac{dy(z)}{dz}+(z^2-\nu^2)y(z)=z^{\mu+1}
\end{equation}
determined by the behavior $y(z)\sim \frac{z^{\mu+1}}{(1+\mu)^2-\nu^2}(1+O(z^2))$ around $z=0$ when $(\mu+n)^2\neq \nu^2$, $n\in \mathbb{N}$. One has \cite{W}
\begin{equation}\label{Lom}
s_{\mu,\nu}(z)=\frac{z^{\mu+1}}{(\mu+1)^2-\nu^2} {}_1F_2\left(1; \frac{\mu-\nu+3}{2},\frac{\mu+\nu+3}{2};-\frac{1}{4}z^2\right).
\end{equation}
The values of $(\mu,\nu)$ corresponding to the parameters $(\frac{4}{3}, \frac{5}{3})$ in (\ref{dz0}) are $(\mu,\nu)=(0,1/3)$. Further, $s_{0,v}(z)$ possesses a nice integral representation, given by
 \begin{equation}\label{int}
s_{0,\nu}(z)=\frac{1}{1+\cos(\pi\nu)}\int_0^\pi \sin(z\sin(t))\cos(\nu t)dt.
\end{equation}
 Indeed, by inserting into the left hand side of (\ref{nhB}) the integral (\ref{int}) and multiplying by $(1+\cos(\pi\nu))$ one gets
\begin{equation}
\begin{split}
& (1+\cos(\pi\nu))\left(z^2\frac{d^2y(z)}{dz^2}+z\frac{dy(z)}{dz}+(z^2-\nu^2)y(z)\right)=\\
&=\int_0^\pi \cos(\nu t)\left(\sin(z\sin(t))(\cos(t)^2z^2-\nu^2)+z\cos(z\sin(t))\sin(t)\right)dt=\\
&=-\left (\cos(\nu t)\cos(z\sin(t))\cos(t)z+\sin(z\sin(t))\sin(\nu t)\nu)\right|_0^\pi=z(1+\cos(\pi\nu)).
\end{split}
\end{equation}
The integral on the right hand side of (\ref{int}) has also the right behavior around $z=0$, i.e. $\sim \frac{z}{1-\nu^2}(1+O(z^2))$, showing the equivalence stated in (\ref{int}). From this integral representation it follows another formula, that will be important in the following. The formula is
\begin{equation}\label{int1}
s_{0,\nu}(z)=\frac{1}{1+\cos(\pi\nu)}\int_0^1 \sin(z t)\frac{\cos(\nu\arcsin(t))+\cos(\nu\pi-\nu\arcsin(t))}{\sqrt{1-t^2}}dt,
\end{equation}
where the range of $\arcsin$ is $(-\pi/2,\pi/2)$. At this point it is useful to recall a Theorem of Polya about the zeros of entire functions defined by integrals like (\ref{int1}). One has
\begin{thm}
[Polya, 1918 \cite{P}] Suppose that the function $f(t)$ is positive and not decreasing in $(0,1)$. Then the function of $z$ defined by
\begin{equation}\label{P}
V(z)=\int_0^1\sin(zt)f(t)dt,
\end{equation}
possesses only real zeros. Further,  if $f(t)$ grows steadily these zeros are simple and the intervals $(k\pi, (k+1)\pi)$, $k=1, 2, \ldots$ contain the positive zeros of $V(z)$, each interval containing just one zero.
\end{thm}
By growing steadily it is meant that the function is not piecewise constant with a finite number
of rational points of discontinuity in $(0,1)$. We can apply the theorem of Pyola to formula (\ref{int1}). It is possible to show indeed that the function multiplying $\sin(zt)$ in (\ref{int1}) is positive and increasing for $|\nu|<1$ and $t\in (0,1)$. Putting together equations (\ref{dz0}), (\ref{Lom}) e (\ref{int1}) we get the following 
\begin{propn}\label{prop23}
The function $\tau_d(z,0)$ defined by the series (\ref{dz0}) has a double zero in $z=0$. All other zeros are simple and are located on the three rays $(e^{\textrm{i}\pi}, e^{\pm\textrm{i}\pi/3})$. The modulus of the simple zeros are contained in the intervals 
\begin{equation}
\left[\left(\frac{3}{2}\pi k\right)^{\frac{2}{3}},\left(\frac{3}{2}\pi (k+1)\right)^{\frac{2}{3}}\right], \quad k=1, 2, \ldots 
\end{equation}
\end{propn}
This global result can be compared with the asymptotic distribution of the zeros. Indeed, since ${}_1F_2\left(1; \frac{3-\nu}{2},\frac{3+\nu}{2}; -\frac{z^2}{4} \right)\sim \frac{C_{\nu}}{z^{3/2}}\cos(z-\frac{3\pi}{4})$ \cite{LW}, where $C_{\nu}$ is a constant, it follows that the modulus of the zeros is asymptotically approximated by 
\begin{equation}\label{za}
|z_k|\sim \left(\frac{3}{2}\pi \left(k+\frac{1}{4}\right)\right)^{\frac{2}{3}},\quad k=1, 2, \ldots
\end{equation}
Formula (\ref{za}) gives quite accurate results also for small values of $k$: indeed If $z_{app,k}$ and $z_{e,k}$ are the approximated (through formula (\ref{za})) and exact values for the modulus of the zeros of $\tau_d(z,0)$, then numerically we find that the relative errors are bounded by the following inequalities: 
\begin{equation}
|1-\frac{z_{app,k}}{z_{e,k}}|< 0.01, \quad \textrm{for}\; k>3.
\end{equation}
\begin{rem}
It can be interesting to remark that the Polya's Theorem and equation (\ref{int1}) imply that the generalized hypergeometric function ${}_1F_2\left(1; \frac{3-\nu}{2},\frac{3+\nu}{2}; z \right)$ has only real negative simple zeros for $|\nu|$<1. These zeros are contained in the intervals
\begin{equation}
-\left[\frac{\pi^2(k+1)^2}{4}, \frac{\pi^2k^2}{4}\right], \quad k=1, 2, \ldots
\end{equation} 
\end{rem}

\section{The principal family}
So far, we characterized the zeros of the particular solution $\tau_d(z,0)$ of equation (\ref{eq1}) by showing that this solution possesses just one double zero, the other being simple. 
Actually, as it will be shown, this property is shared by all solutions
$\tau_d(z,p)$, for any $p \in \mathbb{C}$. The particular solutions $\tau_d(z,p)$ of equation (\ref{eq1}) are characterized by the presence of a double zero in the arbitrary point $z=p$. As we said, 
the position of one zero can be considered a constant of integration. The other constant is the value of the derivative of the function at this zero. The functions $\tau_d(z,p)$ however have a derivative
equal to 0 in $z=p$ (for any $p$). The set of functions $\tau_d(z,p)$ then represents a family of particular solutions of equation (\ref{eq1}). This point will be further clarified in this section. Besides the family $\tau_d(z,p)$, equation (\ref{eq1}) possesses a \emph{principal family} of solutions. Proposition (\ref{p1}) can be easily extended to the case of a solution with a simple zero in an arbitrary point $z=q$. Indeed one has
\begin{propn}\label{psol}
Given $q \in \mathbb{C}$, equation (\ref{eq1}) possesses a solution $\tau(z,q,\alpha)$ with a simple zero in $z=q$. This solution is defined by the series
\begin{equation}\label{d}
\tau(z,q,\alpha)=\sum_{n=1}f_n(z-q)^n, \quad f_1=\alpha,\; f_2=\frac{c}{2}, \; f_3=\frac{q\alpha}{6},\; f_{n+2}=\frac{qf_n+f_{n-1}}{(n+1)(n+2)}
\end{equation}
For $q$ real, $\tau(z,q,\alpha)$ can also be represented as
\begin{equation}\begin{split}
&\tau(z,q,\alpha)=\pi \alpha\left(\textnormal{Bi}(z)\textnormal{Ai}(q)-\textnormal{Ai}(z)\textnormal{Bi}(q)\right)+\\
&+ c\pi \left(\textnormal{Bi}(z)\int_q^z \textnormal{Ai}(x)dx-\textnormal{Ai}(z)\int_{q}^z\textnormal{Bi}(x)dx\right).
\end{split}\end{equation}
\end{propn}
In the previous formulae, $q$ and $\alpha$ are arbitrary parameters and so the series (\ref{d}), that converges everywhere in the complex plane, can be considered the general solution of equation (\ref{eq1}). At this point, a note of caution may be worth to understand the characteristics of this representation of the general solution and its potential intersections with the representation (\ref{dz}). It is clear that if $\alpha=0$,  the series (\ref{d}) represents a solution with a double zero in $z=q$, like the series (\ref{dz}) with $p=q$. We recall however that each element of the set of solutions (\ref{dz}) $\tau_d(z,p)$ possesses, further to the double zero in $z=p$, only simple zeros. Let $\{\Omega_k(p)\}_{k=1,2,...}$ be the set of these simple zeros. If $q=\Omega_n(p)$ for some $n$ and $p$ and $\alpha=\tau'_d(\Omega_n(p))$, from the existence and uniqueness theorem it follows that the two series (\ref{dz}) and (\ref{d}) coincide. In all other cases, the function $\tau(z,q,\alpha)$ defined by the series (\ref{d}), must have only simple zeros since it cannot be represented by any element of the family $\tau_d(z,p)$.

 The fact that the set of solutions of (\ref{eq1}) possessing only simple zeros is not empty follows from the observation that, if it would be empty, the set of particular solutions (\ref{dz}) would represent the general solution of equation (\ref{eq1}), that is impossible since in (\ref{dz}) there is just one free parameter. So we arrived at the conclusion that equation (\ref{eq1}) possesses two family of solutions: the particular family, given by the solutions having just one double zero in $z=p$, all other zeros being simple, and the principal family, given by the solutions of equation (\ref{eq1}) having only simple zeros.  This result is summarized in the following
\begin{propn}\label{prop32}
Let $p\in \mathbb{C}$ be arbitrary and let $\tau_d(z,p)$ be the generic element of the set of particular solutions of equation (\ref{eq1}) defined by the series (\ref{dz}). Let $\{\Omega_k(p)\}_{k=1,2,...}$ be the set of the zeros of $\tau_d(z,p)$. Then, if $q\neq \Omega_n(p)$ for some $n$ and $p$ and $\alpha \neq \tau'_d(\Omega_k(p))$, the series (\ref{d}) gives a generic element of the principal family of solutions of equation (\ref{eq1}). Further, any element of the principal family possesses only simple zeros.
\end{propn}
\begin{rem}
In the case $p$ and $q$ are real, then from proposition (\ref{psol}) it follows that the function defined by the series (\ref{d}) possesses a double zero in $z=p$ if and only if $\alpha$ is the common value of the following expressions:
\begin{equation}
\alpha=\frac{c\int_q^p\textnormal{Ai}(x)dx}{\textnormal{Ai}(q)}=\frac{c\int_q^p\textnormal{Bi}(x)dx}{\textnormal{Bi}(q)}.
\end{equation}
\end{rem}

The characterization of the two families of solutions can be further clarified by looking at the logarithmic derivative of $y(z)$, i.e. to $u(z)=\frac{y'(z)}{y(z)}$. This function is meromorphic and, from (\ref{eq1}), it solves the following differential equation
\begin{equation}\label{equ}
\frac{d^2u(z)}{dz^2}+3u(z)\frac{du(z)}{dz}+u(z)^3=1+zu(z)
\end{equation}
Clearly, $u(z)$ possesses the Painlev\'e property, but still it is instructive to apply the Painlev\'e test to $u(z)$. Here we will use the standard terms that can be found in modern literature about the Painlev\'e test (see e.g. \cite{Hone}). The dominant balances gives singularities of the type $c_0(z-p)^{-1}$, but, seeking the resonances, one discovers that there are two families of solutions:
\begin{itemize}
\item One family characterized by $c_0=1$ and with resonance polynomial given by $(r-1)(r+1)$.
\item The other characterized by $c_0=2$ and with resonance polynomial given by $(r+1)(r+2)$.
\end{itemize}
The two resonances $r=-1$ corresponds to the arbitrariness of the position of the pole for $u(z)$, i.e. the arbitrarines of the position of the zero $z=p$ for $y(z)$. For the other resonances one has:
\begin{itemize}
\item In the first family, $c_0=1$ implies that the zero in $z=p$ of $y(z)$ is simple, whereas the resonance $r=1$ indicates that the coefficient of $(z-p)^0$ in the Laurent expansion of $u(z)$ is the other arbitrary constant describing the solutions of the second order equation (\ref{equ}). This is the principal family of solutions \cite{CFP}.
\item In the second family, $c_0=2$ implies that the zero in $z=p$ of $y(z)$ is double, whereas the resonance $r=-1$ is negative: this indicates that the second family is a particular solution of equation (\ref{equ}) \cite{CFP}.
\end{itemize} 
As we shall see the coefficients of the Laurent expansion of the logarithmic derivative of $y(z)$ are explicitly connected with the zeros of $y(z)$: this connection was crucial for the characterization of the distribution of  the zeros of the solutions of the homogeneous Airy equation in \cite{Z1} \cite{Z2}. The same line of reasoning will be adopted in the next section for the nonhomogeneous case. Before, as will be explained, we need to slightly modify equation (\ref{eq1}) and insert two more parameters.

\section{A mapping among solutions.}\label{sec4}
 As we have shown in the previous sections, it is possible to pick a solution of equation (\ref{eq1}) by assigning the
position of one zero and the value of the first derivative of the function at this zero. If the solution of equation (\ref{eq1}) has a zero at $z=z_0$,
with $z_0$ an arbitrary point, by a shift $z\to z+z_0$ and a rescaling of $z$, we get the following equation for $y(z)$:
\begin{equation}\label{eq2}
\frac{d^2y(z)}{dz^2}=(az+b)y(z)+c.
\end{equation}
A solution of equation (\ref{eq2}) with a zero in $z=0$ corresponds to a solution of equation (\ref{eq1}) with a zero in  $z_0=\frac{b}{\eta^2}$, where $a=\eta^3$. Instead of looking at solutions 
of equation (\ref{eq1}) with zeros in $z=z_0$, with $z_0$ arbitrary, we can look at the particular solution of equation (\ref{eq2}) with a zero in $z=0$ for an arbitrary value of the parameter $b$. This second point
of view will be adopted in the rest of the paper. 

As we have seen, the particular solutions $\tau_d(z,0)$ (\ref{dz}) of equation (\ref{eq1}) possesses just one double zero in $z=0$, with all other zeros simple. Equivalently, for $b=0$ there is a solution of equation (\ref{eq2}) possessing a double zero in $z=0$, the other zeros being simple. We would like to prove that this property is shared by all the particular solutions $\tau_d(z,p)$ (\ref{dz}) of equation (\ref{eq1}). This is equivalent to say that, for any $b\neq 0$, equation (\ref{eq2}) possesses a solution with a double zero in $z=0$, with all the other zeros simple. We can give two arguments to justify the presence of just one double zero. For clarity, let us call the solution of equation (\ref{eq2}) with a double zero in $z=0$ $c\Xi(z,a,b)$, where $\Xi(z,a,b)$ is defined by the series
\begin{equation}\label{xi}
\Xi(z,a,b)=\sum_{n=2}d_nz^n, \quad d_2=\frac{1}{2}, \; d_3=0, \; d_4=\frac{b}{24},\; d_{n+2}=\frac{bd_n+ad_{n-1}}{(n+1)(n+2)}.
\end{equation}
By multiplying equation (\ref{eq2}) by $y'(z)$, integrating and taking into account that $\Xi(0,a,b)=0$ and $\Xi'(0,a,b)=0$, we find, for $z \in \mathbb{R}$
\begin{equation}\label{firarg}
\left(\Xi'(z,a,b)\right)^2=(az+b)\Xi^2(z,a,b)+2c\Xi(z,a,b)-a\int_0^z\Xi^2(z,a,b)dz
\end{equation}
In the case we have another double zero in $z=p$, $p\in \mathbb{R}$, one should have
\begin{equation}\label{imp}
\int_0^p\Xi^2(z,a,b)dz=0,
\end{equation}
which is impossible. Notice that we tacitly assumed $\Xi$ to be real on the real axis and this is the case if the coefficients appearing in equation (\ref{eq2}) are real.

For the second argument, we observe that the solutions of equation (\ref{eq2}) with a double zero in $z=0$ have a nice homogeneity property with respect to the parameters. Indeed, for any $\lambda \in \mathbb{C}^{*}$ one has:
\begin{equation}\label{hity}
\Xi(\lambda^{-1}z,\lambda^3a,\lambda^2b)=\lambda^2\Xi(z,a,b).
\end{equation}
For $b=0$, from (\ref{dz0}) we have
\begin{equation}\label{dz0d}
\Xi(z,a,0)=\sum_{n=0}\frac{(3a)^n n!}{(3n+2)!}z^{3n+2} = \frac{z^2}{2} {}_1F_2\left(1; \frac{4}{3},\frac{5}{3};\frac{az^3}{9}\right).
\end{equation}
Now we ask how the zeros of $\Xi(z,a,0)$ will move in the complex plane if we take a small value of $b$. We remember that the zeros of $\Xi(z,a,0)$, for $a=1$ are distributed on the rays $(e^{\textrm{i}\pi}, e^{\pm\textrm{i}\pi/3})$. Changing this value of $a$ corresponds to a scaling (for $|a| \neq 1$) and global rotation (for arg$(a) \neq 0$) of their positions. 
If we  expand $\Xi(z,a,b)$ in a series of $b$, each coefficient being a function of $z$ and $a$, we get
\begin{equation}\label{xis}
\Xi(z,a,b)=\Xi_0(z,a)+b\Xi_1(z,a)+O(b^2),
\end{equation}
 where $\Xi_0(z,a)=\Xi(z,a,0)$ is given by (\ref{dz0d}). The function $\Xi_1(z,a)$ solves the differential equation
 \begin{equation}
 \frac{d^2\Xi_1(z,a)}{dz^2}=az\Xi_1(z,a)+\Xi_0(z,a),
 \end{equation}
with the initial conditions $\Xi_1(0,a)=0$, $\Xi'_1(0,a)=0$. It is possible to show that the corresponding solution is defined by the series
\begin{equation}\label{d1}
\Xi_1(z,a)=\sum_{n=0}a^n\left(\frac{3^{n+1} (n+1)!}{(3n+4)!}-\frac{1}{3^{n+1}n!}\prod_{k=1}^{n+1}\frac{1}{3k+1}\right)z^{3n+4},
\end{equation}
that can be resummed to the difference of two generalized hypergeometric functions:
\begin{equation}\label{d2}
\Xi_1(z,a)= \frac{z^4}{8} {}_1F_2\left(1; \frac{5}{3},\frac{7}{3};\frac{az^3}{9}\right)-\frac{z^4}{12}{}_1F_2\left(1; 2,\frac{7}{3};\frac{az^3}{9}\right) .
\end{equation}
If $\xi_{1,k}$ are the zeros of $\Xi(z,a,b)$ and $\xi_{k,0}$ are those of $\Xi(z,a,0)$, at first order in $b$  we have
\begin{equation}\label{xir}
\xi_{1,k}-\xi_{0,k}=-b\frac{\Xi_1(z_{0,k},a)}{\Xi'_0(z_{0,k},a)},
\end{equation}
or, in terms of hypergeometric functions:
\begin{equation}
\xi_{1,k}-\xi_{0,k}=-\frac{5b}{9a}\frac{3\, {}_1F_2\left(1; \frac{5}{3},\frac{7}{3};\frac{az_{0,k}^3}{9}\right)-2 \,{}_1F_2\left(1; 2,\frac{7}{3};\frac{az_{0,k}^3}{9}\right)}{{}_1F_2\left(2; \frac{7}{3},\frac{8}{3};\frac{az_{0,k}^3}{9}\right)},
\end{equation}
Notice that $\xi_{1,k}-\xi_{0,k}$ depends only on the value of $\xi_{0,k}^3$: the difference is the same for all the triple of zeros lying on the circle of radius $|\xi_{0,k}|$. All the zeros move in the same direction and, further, each triple of zeros on a given circle move of the same amount in that direction, like if they belong to a rigid body. The denominator of (\ref{xir}) is always different from zero since, from proposition (\ref{prop23}), the zeros of $\Xi_0$ are simple. The numerator is an entire function of $\xi_{0,k}$. It follows that, for any $k$, the ratio $\frac{\Xi_1(z_{0,k},a)}{\Xi'_0(z_{0,k},a)}$ is bounded by a suitable constant $M_k$. Then, for any $k$, it is possible to choose a constant $b_k$ such that, for $|b|<|b_k|$, the zeros $\xi_{1,n}$, $n<k$, are all simple. This implies that, for suitable small values of $b$, the zeros of  $\Xi(z,a,b)$ are simple (apart the double zero in $z=0$). This line of reasoning is independent on the value of $a$, since, as we said, the zeros of $\Xi_0(z,a,0)$ are only scaled and/or rotated with respect to those of $\Xi(z,1,0)$. Now we can go back to the homogeneity property (\ref{hity}). Indeed, we have seen that the zeros of $\Xi(z,\frac{a}{\lambda^3}, 0)$ are simple for any value of the ratio $a/\lambda^3$. Also, the zeros of  $\Xi(z,\frac{a}{\lambda^3}, b)$ are simple for sufficiently small values of $b$. From the homogeneity property (\ref{hity}) it follows that the zeros of $\Xi(z, \frac{a}{\lambda^3},b)$ coincide with those of  $\Xi(\lambda^{-1} z, a, \lambda^2 b)$ for any value of $\lambda$, i.e. coincide, with a proper rescaling and rotation, with the zeros of  $\Xi(z, a, \lambda^2 b)$. Since $\lambda$ is arbitrary we get that the zeros of $\Xi(z,a,b)$ are all simple (apart $z=0$) for any value of $a$ and $b$. 

Up to now we investigated the distribution of the zeros of the solutions of equation (\ref{eq1}) with a double zero in $z=p$. In the rest of the paper we are going to investigate the distribution of the zeros of the solutions of equation (\ref{eq1}) possessing only simple zeros (equivalently, the solution of equation (\ref{eq2}) for suitable values of the parameters $a$ and $b$). To distinguish between the two types of solutions, we make the following definition:
\begin{defn} $X(z,a,b,c)$ are the solutions of (\ref{eq2}) having a simple zero in an arbitrary point of the complex plane, all the other zeros being simple. These zeros are denoted by $\chi_k(a,b,c)$. c$\Xi(z,a,b)$ are the solutions of (\ref{eq2}) having one double zero in an arbitrary point of the complex plane. The other zeros are denoted by $\xi_k(a,b)$. 
\end{defn}
It will be useful the following remark:
\begin{rem}\label{gs}
If $y_0(z,a,b,c)$ is any particular solution of equation (\ref{eq2}), then another solution is given by T$_{A,B}y_0$, where T$_{A,B}$ is defined by
\begin{equation}\label{trans}
\mathrm{T}_{A,B}y_0(z,a,b,c)=Ay_0(z-B,a,b+aB,\frac{c}{A})
\end{equation}
with $A$ and $B$ two arbitrary constants.
\end{rem}
 From proposition (\ref{prop32}) and the discussion just before that proposition, it is possible to write the following charts for the action of the operator T$_{A,B}$ on the two family of solutions $X(z,a,b,c)$ and $\Xi(z,a,b)$:
 \begin{figure}[h]
%\caption{A picture of a gull.}
  \centering
    \includegraphics[scale=0.5]{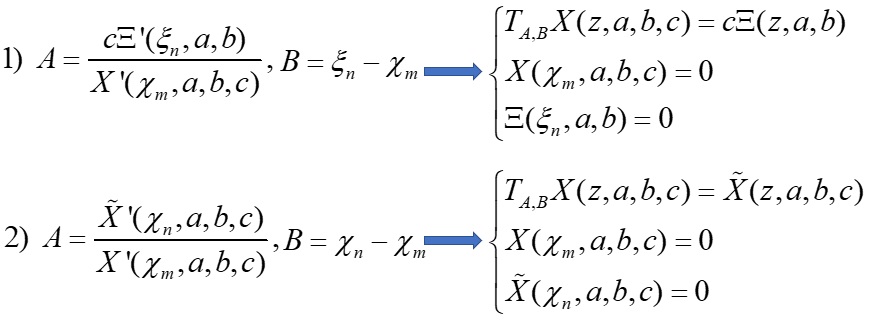}
\end{figure}
The transformation 1) brings the element of the family $X(z,a,b,c)$ having a zero in $z=\xi_m$ to an element of the family $\Xi(z,a,b)$ having a zero in $z=\xi_n$: if this zero is simple, then $A\neq 0$ and the function $\Xi(z,a,b)$ has a double zero in some other point of the complex plane. If $\xi_n$ is the double zero, then $A=0$. The transformation $T_{0,B}$ eliminates all the terms not proportional to $c$ in $X(z,a,b,c)$. We recall that the element of the family $X(z,a,b,c)$ having a zero in an arbitrary point $z=\chi_m$ is linear in $c$ and in the first derivative of $X$ evaluated in $\chi_m$, as can be seem also from proposition \ref{psol})
 
 The transformation 2) brings the element of the family $X(z,a,b,c)$ having a zero in $z=\xi_m$ to an element $\tilde{X}(z,a,b)$ of the same family. This new solution has a simple zero in $z=\chi_n$. 
 
 It is also possible to get a transformation from $\Xi(z,a,b,c)$ to the same family: the value of $A$ is redundant in this case since the starting solution is proportional to $c$ and there is just the parameter $B$ surviving in equation (\ref{trans}).
 
 In the next it will be useful to answer to the following question: it is possible, for a suitable choice of $A$ and $B$ in the transformation 2), to get the identity transformation? If we start from the solution possessing a zero in $z=\xi_m$, then the value of $B$ must be equal to $\xi_n-\xi_m$, where $\xi_n$ is another zero of the \emph{same} solution. Also, the value of $A$ must be equal to the ratio of the first derivative evaluated at $z=\xi_n$ over the first derivative evaluated at $z=\xi_m$. Since $\xi_n$ is an arbitrary zero, actually there is an infinite, numerable set of choices for the values of $A$ and $B$. In this case $T_{A,B}$ is the identity transformation and we can write the following
\begin{propn}\label{propu}
If $X(z,a,b,c)$ is the solution of equation (\ref{eq2}) having a zero in $z=\xi_m$, then one has
\begin{equation}
X(z,a,b,c)=A_{n,m}X\left(z-(\xi_n-\xi_m),a,b+a(\xi_n-\xi_m),c/A_{n,m}\right)
\end{equation}
where $A_{n,m}=\frac{X'(\xi_n,a,b,c)}{X'(\xi_m,a,b,c)}$ and $\xi_n$ is any other zero of $X(z,a,b,c)$.
\end{propn}
 As we have seen, a solution of equation (\ref{eq2}) with a zero in $z=0$ corresponds to a solution of equation (\ref{eq1}) with a zero in $z_0=\frac{b}{\eta^2}$, where $a=\eta^3$. Also, the characterization of the solutions of equation (\ref{eq1}) with a zero in an arbitrary point $z_0$ of the complex plane corresponds to the characterization of the particular solution of equation (\ref{eq2}) having a zero in $z=0$ for arbitrary values of $a$ and $b$. The proposition (\ref{propu}) gives  a quasi-periodic property of the zeros of this particular solution. Indeed, from (\ref{propu}) we get the following
\begin{rem}\label{r4}
Let $S(z,a,b,c)$ be a particular solution of equation (\ref{eq2}) having a simple zero in $z=0$, with all the other zeros simple. Let $\{\xi^0_k(a,b,c)\}$ denote the set of the zeros of $S(z,a,b,c)$. Then from proposition (\ref{propu}) it follows that, for any choice of $\xi^0_n$ one has
\begin{equation}\label{z}
\left\{\xi^0_k\left(a,b+a\xi^0_n(a,b,c),c/A_n\right)\right\}=\left\{\xi^0_k(a,b,c)-\xi^0_n(a,b,c)\right\}
\end{equation}
where $A_n=\frac{S'(\xi_n,a,b,c)}{S'(0,a,b,c)}$.
\end{rem}
We notice that for $c=0$ formula (\ref{z}) gives the quasi-periodicity of the zeros of the solutions of the homogeneous Airy equation, as described in \cite{Z1} \cite{Z2}, whereas for $c=0$ and $a=0$ one gets the periodicity of the trigonometric functions.

 \section{A cubic recursion for the zeros}\label{sec5}
To further characterize the zeros of the function $S(z,a,b,c)$, we can look at the logarithmic derivative of this function. Indeed, since $S(z,a,b,c)$ is an entire function with order of growth equal to $3/2$, from the Weierstrass-Hadamard factorization theorem \cite{T} we get
\begin{equation}\label{WH}
 S(z,a,b,c)=S'(0,a,b,c)ze^{\beta z}\prod_{k=1}\left(1-\frac{z}{\xi_k}\right)e^{\frac{z}{\xi_k}} 
\end{equation}
where $\beta=\frac{c}{2S'(0,a,b,c)}$ and hereafter we omit the apex ``0'' from the zeros $\xi^0_n$ for ease of readability. The logarithmic derivative of $S$, $u=\frac{S'}{S}$, is a meromorphic function and, from the previous product fromula, possesses the (global) Mittag Leffler representation
\begin{equation}\label{ML}
u(z,a,b,c)=\frac{1}{z}+\beta+\sum_{\xi_n\neq 0}\frac{1}{z-\xi_n}+\frac{1}{\xi_n}
\end{equation}
From the differential equation (\ref{eq2}) it is possible to show that $u$ solves the following second order, cubic, differential equation
\begin{equation}\label{equ}
\frac{d^2u(z)}{dz^2}+3u(z)\frac{du(z)}{dz}+u(z)^3=a+(az+b)u(z).
\end{equation}
We can expand $u(z,a,b,c)$ in a (local) Laurent series around $z=0$, getting
\begin{equation}\label{L}
u(z,a,b,c)=\frac{1}{z}+\sum_{n=0}c_nz^n.
\end{equation}
From the differential equation (\ref{equ}), the coefficients $c_n(a,b,c)$ solve the following recursion
\begin{equation}\label{cn}\begin{split}
& (n+4)(n+2)c_{n+2}=bc_n+ac_{n-1}-3\sum_{k=0}^{n+1}c_kc_{n+1-k}+\\
&-3\sum_{k=0}^{n}(k+1)c_{k+1}c_{n-k}-\sum_{k=0}^{n}\sum_{j=0}^kc_jc_{k-j}c_{n-k}, \quad n\geq 1
\end{split}\end{equation}
%\endgroup
where $c_0=\beta, c_1=\frac{b}{3}-c_0^2, c_2=c_0^3-\frac{1}{4}bc_0+\frac{a}{4}.$ The compatibility between the expansion (\ref{L}) and the representation (\ref{ML}) gives the following expression relating the coefficients $c_n$ to the zeros $\xi_k$:
\begin{equation}\label{thecn}
c_n=-\sum_{\xi_k\neq 0}\frac{1}{\xi_k^{n+1}}, \qquad n=1,2...
\end{equation}
In the case $c=0$ (i.e. for the solutions of the homogeneous Airy equation) the previous recursion is not more cubic but quadratic in the coefficients $c_n$. In that case, formula (\ref{cn}) can be inverted, in the sense that known all the coefficients $c_n$ it is possible, in principle, to calculate all the zeros $\xi_k$ \cite{Z1} \cite{Z2}. In this case we can repeat the same arguments as given in \cite{Z1} and we refer to that work for more details. The result is the following
\begin{propn}\label{prrec}
Let $c_n(a,b+ax,\beta)$ be the sequence of polynomials defined by the recursion (\ref{cn}) with $b\to b+ax$.
%\begingroup\small
%\begin{equation}\begin{split}
%& (n+4)(n+2)c_{n+2}=(b+ax)c_n+ac_{n-1}-3\sum_{k=0}^{n+1}c_kc_{n+1-k}+\\
%&-3\sum_{k=0}^{n}(k+1)c_{k+1}c_{n-k}-\sum_{k=0}^{n}\sum_{j=0}^kc_jc_{k-j}c_{n-k}, \quad n\geq 1
%\end{split}\end{equation}
%\endgroup
%where $c_0=\beta, c_1=\frac{b+ax}{3}-c_0^2, c_2=c_0^3-\frac{1}{4}(b+ax)c_0+\frac{a}{4}.$
Then, if the distance between successive zeros is decreasing and $\xi_1$ is the zero of the function $S(z,a,b,c)$ closest to the origin, a subset of the zeros of $S(z,a,b,c)$ are given  
by
\begin{equation}\label{rec}
\xi_{k+1}=\xi_k + \lim_{n\to \infty} \frac{c_{n}(a,b+a\xi_k,\beta/S'(\xi_k))}{c_{n+1}(a,b+a\xi_k,\beta/S'(\xi_k))}, \quad \xi_0=0
\end{equation}
\end{propn}
The subset of the zeros are those on the semi-axis containing both $0$ and $\xi_1$ in the direction of $\xi_1$. We are assuming also that there is just one zero closest to the origin. On the contrary, with more than one zero closest to the origin lying on the circle of radius $|\xi_1|$, the right hand side of (\ref{rec}) would oscillate indefinitely. The constraint that the distance between successive zeros must be decreasing  is crucial: for $c=0$ (i.e. $\beta=0$ in the recursion (\ref{rec})) it has been shown \cite{Z1} that, for suitable choice of the constants $a$ and $b$, this property holds true and the recursion (\ref{rec}) can be effectively used to calculate the corresponding zeros. In the case $c\neq 0$ the situation is different: for small values of $c$ we expect that the this property may hold true again, but, by increasing $c$, we have a different picture. Also, another zero can appear on the  semi-axis opposite to that containing $0$ and $\xi_1$. This rich and interesting portrait will be further considered in a next work. Here we will give just few examples of application of the recursion (\ref{rec}). We know from \cite{Z1}-\cite{Z2} that, for $c=0$, $a<0$ and $b<0$ $S(z,a,b,0)$ has an infinite number of zeros on the positive real axis. We first consider small values of the parameter $c$: for definiteness, we take $(a,b,c)=(-1,-1,-0.1)$. After, we will consider a value of $c$ comparable to those of $a$ and $b$: we will take  $(a,b,c)=(-1,-1,-1)$. In both the examples we fix the value of $S'(0,a,b,c)$ to be equal to $1$ (the distribution of the zeros depends only on the ratio between $c$ and $S'(0,a,b,c)$). In the last case, the zero closest to $z=0$ is real and is equal to 
\begin{equation}
\xi_1=1.4230603\ldots
\end{equation}
The next real zero is about equal to $\xi_2=3.53175\ldots$, meaning that the condition given in the Proposition (\ref{prrec}) about the distance of successive zeros is not satisfied, since $|\xi_2-\xi_1|>|\xi_1-\xi_0|$: \emph{in this case the recursion (\ref{rec}) gives the zero that is closest to $\xi_1$}, i.e. one gets
\begin{equation}
\xi_{2}=0=\xi_1 + \lim_{n\to \infty} \frac{c_{n}(-1,-1-\xi_1,-1/2)}{c_{n+1}(-1,-1-\xi_1,-1/2)},
\end{equation} 
that is a 2-peiodic sequence. On the contrary, in the first case, i.e. for $(a,b,c)=(-1,-1,-0.1)$,  the conditions of Proposition (\ref{prrec}) are satisfied: the zero closest to $\xi_0=0$ is
\begin{equation}
\xi_1=2.0977152\ldots
\end{equation}
The next two zeros, found with the recursion (\ref{rec}) itself, are given by
\begin{equation}
\xi_2=3.7233151\ldots,\quad \xi_3=5.0507149\ldots
\end{equation}
The previous zeros have been obtained by taking the ratio of $c_n/c_{n+1}$ up to $n=80$: with larger number of $n$ it would be possible to get more significant digits in the calculation of the zeros.

\section{Conclusions}
We hope to have highlighted some of the very interesting mathematical structures defined by the general solutions of equation (\ref{eq1}) in the complex plane. The possibility to look at a particular solution of equation (\ref{eq2}) to characterize the general solution of equation (\ref{eq1}) is surely useful and may be relevant in future potential developments of this work. The dynamics of the zeros by varying the parameters in equation (\ref{eq2}), in the spirit of \cite{Z2} , is an issue that will be investigated in future works. The relation between nonlinear recurrence relations, like (\ref{cn}), and the zeros of entire functions, seems to be another point deserving more consideration in the literature.

\section*{Acknowledgments} The support of University of Brescia, INdAM-GNFM and INFN is gratefully acknowledged. Also, the author wishes to thank Prof. O. Ragnisco for his advice and comments.

\end{document}